\documentclass[12pt,reqno,oneside]{amsart}
\usepackage{amsmath,amsthm,amsfonts,amssymb}
\usepackage[mathscr]{eucal}
\usepackage{indentfirst}
\usepackage{url}

\theoremstyle{plain}
\newtheorem{teo}{Theorem}[section]
\newtheorem{cor}[teo]{Corollary}
\newtheorem{lem}[teo]{Lemma}

\theoremstyle{definition}
\newtheorem{defn}[teo]{Definition}
\newtheorem{exa}[teo]{Example}

\numberwithin{equation}{section}

\def\bbP{{\mathbb P}}
\def\bbZ{{\mathbb Z}}
\def\bbN{{\mathbb N}}

\def\bbT{{\mathbb T}}

\def\qed{\hfill $\square$}

\def\cpp{\textit{Cone Percolation Model}}

\begin{document}
\baselineskip=22pt
\title[The cone percolation model]{The cone percolation model on Galton-Watson and on spherically symmetric trees}
\author{Valdivino~V.~Junior}
\author{F\'abio~P.~Machado}
\author{Krishnamurthi~Ravishankar}
\address[F\'abio~P.~Machado]
{Institute of Mathematics and Statistics
\\ University of S\~ao Paulo \\ Rua do Mat\~ao 1010, CEP
05508-090, S\~ao Paulo, SP, Brazil - fmachado@ime.usp.br }

\noindent
\address[Valdivino~V.~Junior]
{Federal University of Goias
\\ Campus Samambaia, Goi\^ania, GO, Brazil - vvjunior@ufg.br}

\address[Krishnamurthi Ravishankar]
{NYU-ECNU Institute of Mathematical Sciences at NYU Shanghai, 3663
Zhongshan Road North, Shanghai 200062 and 1555 Century Ave, Pudong, Shanghai, China - kr26@nyu.edu}

\noindent

\thanks{Research supported by CNPq (310829/2014-3), FAPESP (09/52379-8), PNPD-Capes 536114 and Simons Foundation Collaboration grant 281207}



\date{\today}

\begin{abstract}
We study a rumour model from a percolation theory and branching process point of view. The 
existence of a giant component is related to the event where the rumour, which started from 
the root of a tree, spreads out through an infinite number of its vertices. We present 
lower and upper bounds for the probability of that event, according to the distribution of the 
random variables that defines the radius of influence of each individual. We work with Galton-Watson branching trees (homogeneous and non-homogeneous) and
spherically symmetric trees which includes homogeneous and $k-$periodic trees.
\end{abstract}

\maketitle

\section{Introduction and basic definitions}
\label{S: Introduction}

Lebensztayn and Rodriguez~\cite{LR}, introduced a disk percolation model on general graphs
where a reaction chain starting from the origin of the graph, based on independent 
copies of a geometric random variables, may lead to the existence of a giant component. 

This line of research was continued by Junior {\it et al}~\cite{Junior} and~\cite{Junior2}, focusing
on ${\mathbb N}$ and on the homogeneous tree respectively, studing a family of dependent long range (not necessarily homogeneous)
percolation model. They studied the criticality of each model, presenting suficient conditions 
under which the processes reach a giant component with positive probability.
Besides they presented bounds 
for the probability of having a giant component based on what they considered the 
radius of influence of each vertex of ${\mathbb N}$.

Gallo {\it et al}~\cite{GGJR} computed precisely the probability of having a giant component
for the homogeneous version of one of the models proposed in Junior {\it et al}~\cite{Junior}, 
and obtained information about the distribution of the range of the cluster of the origin when it is finite.  
Besides that, they obtained a law of large numbers 
and a central limit theorem for the proportion of the cluster of the origin in a range of size $n $
as $n$ diverges. The key step of of proofs presented in Gallo {\it et al}~\cite{GGJR} is to show 
that, in each model, the vertices belonging to the cluster of the origin can be related to a 
suitably chosen discrete renewal process. Related results have been obtained recently by 
Bertachi and Zucca~\cite{BZ}. All these research papers are to a different degree, stimulated by the seminal work of Benjamini and Schram~\cite{BS} when they proposed the study of percolation theory beyond the nearest neighbor independent setup on $\bbZ^d$.

Here we focus on Galton-Watson, homogeneous, periodic and spherically symmetric trees
in a process where the \textit{radius of influences} is given by non-negative discrete random variables. In the paper we use the letter $R$ to refer to that
random variable and to make formulas neater we define
$p_k = \bbP(R=k) \mbox{ for } k=0,1,\dots $
To avoid trivialities we assume throughout this paper that $p_0 \in (0,1).$ 
A graph $G$ is said a \emph{tree} if for
any pair of its vertices there is one and only one path (a subset of edges) conecting them.
By $|A|$ we denote the cardinality of $A$.
The \emph{degree} of a vertex is the cardinality of its set of neighbors. For two vertices $u,v$
let $d(u,v)$, be the distance between $u$ and $v$, that is the number of edges
the path from $u$ to $v$ has. 

Consider a tree $\bbT$ (a connected graph with no cycles) and its set of vertices ${\mathcal V}(\bbT)$. 
Single out one vertex from ${\mathcal V}(\bbT)$ and call this ${\mathcal O}$, the origin of ${\mathcal V}(\bbT)$. For each two vertices $u,v \in {\mathcal V}(\bbT)$, 
consider that $ u \leq v$ if $u$ belongs to the path connecting ${\mathcal O}$ to $v$. 

For a tree $\bbT$ and $n \geq 1$ we define

\begin{displaymath}
T^u := \{ v \in \mathcal{V}: u \leq v \},
\end{displaymath}
\begin{displaymath}
T_n^u := \{ v \in T^u: d(v,{\mathcal O}) \leq d(u,{\mathcal O}) + n \}
\end{displaymath}
and
\begin{equation}
\label{E: defPartialT}
M_n(u) := | \partial T_n^u | := |\{ v \in T^u: d(v,{\mathcal O}) = d(u,{\mathcal O}) + n\}|.
\end{equation}

As in Junior {\it et al}~\cite{Junior2}, we say that the process \textit{survives} if the 
number of vertices involved is infinite.
Otherwise we say the process \textit{dies out}. Our main interest is to obtain resuts concerning whether the 
process has positive probability of involving an infinite set of individuals. Besides we present 
lower and upper bounds for the probability of that event, according to the distribution of the 
random variables that defines the radius of influence of each individual. 

The paper is organized as follows. 
Sections~\ref{S: Homogeneous Trees}, \ref{S: Periodic Trees}, 
\ref{S: Spherically Symmetric Trees} 
and \ref{S: Galton-Watson Branching Trees} present the main results and specific setups and distributions 
for the Cone Percolation model on Homogeneous Trees, Periodic Trees, Spherically Symmetric Trees and Galton-Watson T rees respectively.
Section~\ref{S: Proofs} brings the proofs for the main results presented along sections~\ref{S: Homogeneous Trees}, \ref{S: Periodic Trees},
\ref{S: Spherically Symmetric Trees} and \ref{S: Galton-Watson Branching Trees} together with
auxiliary lemmas and useful inequalities.

\section{Homogeneous Trees}
\label{S: Homogeneous Trees}

\noindent
Let us start off with a definition.
\begin{defn} 
\label{D: DefinicaoCPP}
The \cpp\ on $\bbT$.\\
Let $\{R_v\}_{\{ v \in {\mathcal V}(\bbT) \}}$ and $R$
be a set of independent and identically distributed random variables.
Furthermore, for each $u \in {\mathcal V}(\bbT)$, we define the random sets
\begin{equation}
\label{E: defBu}
B_u = \{v \in {\mathcal V}(\bbT): u \leq v \hbox{ and } d(u,v) \leq R_u\}.
\end{equation}
\noindent
With these sets we define the \cpp\ on $\bbT$, the non-decreasing sequence of random 
sets $I_0 \subset I_1 \subset \cdots$
defined as $ I_0 = \{{\mathcal O}\} $ and inductively $I_{n+1} = \bigcup_{u \in I_n} B_u$
for all $ n \geq 0.$ 
\end{defn}

\begin{defn}
The \cpp\ $survival$\\
\label{D: Survival}
Consider $ I = \bigcup_{n \geq 0} I_n$ be the connected component of the
origin of $\bbT$. Under the rumor process interpretation, $I$ is the set of vertices
which heard the rumor. We say that the process \emph{survives} if $|I|=\infty,$
referring to the surviving event as $V.$
\end{defn}

\begin{defn}
Rooted homogeneous trees\\
\label{D: RootedTree}
We say that a tree, $\bbT_d$, is \textit{homogeneous}, if each one of its vertices has degree $d+1$.
From $\bbT_{d}$ we define $\bbT^+_d$, a \textit{rooted homogeneous tree}. Pick a $u \in {\mathcal V}(\bbT_d)$ 
such that $d({\mathcal O},u)=1$ and consider
\[ \bbT^+_d(u) = \{v \in {\mathcal V}(\bbT_d): u \leq v \}.\]
\[ \bbT^+_d := \bbT_d \backslash \bbT^+_d(u) \]
\end{defn}

Consider $\bbP_+$ and $\bbP$ the probability measures associated to the processes on $\bbT_d^+$ and $\bbT_d$
(we do not mention the random variable $R$ unless absolutely necessary). By
a coupling argument one can see that for a fixed distribution of $R$
\begin{equation}
\label{E: NotEqual}
\bbP_+(V) \leq \bbP(V).
\end{equation}

Furthermore, by the definition of $\bbT_d^+$ and its relation with $\bbT_d$ we have that
for a fixed distribution of $R$
\begin{equation}
\label{E: Equal}
\bbP_+(V)=0 \hbox{ if and only if } \bbP(V)=0.
\end{equation}

\begin{teo}
\label{T: THT1}
Consider a \cpp\  on $\bbT_d$.
Then, for $\mathbb{E}(d^R) < 2 - \frac{1}{d}$, we have
\begin{displaymath}
 \frac{d + \mathbb{E}\displaystyle \left (d^R \right) - p_0 }{d[1-\mathbb{E}\displaystyle \left (d^R \right) + p_0]} \leq \mathbb{E}(|I|)  \leq \frac{\mathbb{E}\displaystyle \left (d^R \right)+d-2}{2d - 1 - d\mathbb{E}\displaystyle \left (d^R \right)}.
\end{displaymath}
\end{teo}

\begin{exa} 
Consider $R \sim {\mathcal B}(p)$, a radius of influence sa\-tis\-fying
\[
\bbP(R = 1) = p = 1 - \bbP(R = 0),
\]
with $pd < 1$. Then we have
\begin{displaymath}
\mathbb{E}(|I|) = \frac{1+p}{1-dp}.
\end{displaymath}
\end{exa}

\begin{exa} Consider $R \sim {\mathcal G}(1-p)$, a radius of
influence satisfying
\[
\bbP(R = k) = (1-p)p^{k}, k=0,1,2,\dots
\]
and assume also $pd < \frac{1}{2}$.
So we have
\begin{displaymath}
\frac{1-dp +p -p^2}{1-2dp+dp^2} \leq \mathbb{E}(|I|) \leq \frac{1-dp-p}{1-2dp}.
\end{displaymath}
That gives us a fairly sharp bound even when we pick $p$ and $d$ such that $pd$ is very close to $\frac{1}{2}$ as, for example, $p = 10^{-6}$ and $d = 499,000$. 
For these parameters we get $250.438 \leq \mathbb{E}(|I|) \leq 250.501$.
\end{exa}

\begin{exa}
For $R \sim {\mathcal B}(n,p)$, a radius of influence satisfying
\[
\bbP(R = k) = {n \choose k}p^{k}(1-p)^{n-k}, k=0,1,\dots,n
\]
and $ p < \frac{1}{d-1} [\sqrt[n]{\frac{2d-1}{d}}-1] $, we have
\begin{displaymath}
\frac{d + (dp+1-p)^n - (1-p)^n}{d[1-(dp+1-p)^n + (1-p)^n]} \leq \mathbb{E}(|I|) \leq \frac{(dp+1-p)^n + d - 2}{2d - 1 - d(dp+1-p)^n}.
\end{displaymath}
 Assuming 

$d = 1,000$, $n=2$ and $p = 4 \times 10^{-4}$ we have $24.825 \leq \mathbb{E}(|I|) \leq 24.924.$

\end{exa}

\begin{exa}
For $R \sim {\mathcal P}(\lambda)$, a radius of influence sa\-tis\-fying
\[
\bbP(R = k) = \frac{\exp(-\lambda) \lambda^k}{k!}, k=0, 1, 2, \dots
\] 
and $\lambda < \ln (\sqrt[d-1]{2 - \frac{1}{d}})$, we have
\begin{displaymath}
\frac{d + e^{(d-1)\lambda} - e^{-\lambda}}{d[1 - e^{(d-1)\lambda} + e^{-\lambda}]} \leq \mathbb{E}(|I|) \leq \frac{e^{(d-1)\lambda} + d - 2}{2d - 1 - de^{(d-1)\lambda}}.
\end{displaymath}
In particular, if

$d = 1,000$ and  $\lambda = 6 \times 10^{-4}$, we find $5.613 \le  \mathbb{E}(|I|) \le 5.625$.
\end{exa}

\section{Periodic Trees}
\label{S: Periodic Trees}

\begin{defn}
We define a $k$\textit{-periodic tree} with degree $\tilde{d} = (d_1, \cdots, d_k)$, $d_i \geq 2$ for all 
$i=1,2,\cdots, k$, as tree such that for any vertex whose 
distance to the origin is $nk+i-1$ for some $ n \in\bbN $ has degree $d_i + 1$. 
We refer to this tree as $\bbT_{\tilde{d}}$.
\end{defn}

A few useful quantities to present the results in this section are
\begin{align*}
d_{(i)} &= \textrm {the $i$-th smallest value in } \tilde{d}, \\
G&=G(\tilde{d}) := \sqrt[k]{\prod_{j=1}^{k}}d_{j}, \\
c_0 := 1 \textrm{ and }~c_i :&=
\frac{\prod_{j=1}^{i}d_{(j)}}{\sqrt[k]{\prod_{j=1}^{k}(d_{j})^{i}}}
= \frac{\prod_{j=1}^{i}d_{(j)}}{G^{i}}, i = 1, \cdots, k-1; \\
\bar{c_0} := 1 \textrm{ amd }~ \bar{c_i} :&=
\frac{\prod_{j=k+1-i}^{k}d_{(j)}}{\sqrt[k]{\prod_{j=1}^{k}(d_{j})^{i}}}
= \frac{\prod_{j=k+1-i}^{k}d_{(j)}}{G^{i}}, i = 1, \cdots, k-1.
\end{align*}

\begin{defn}
\label{F: funcaoh}
For $i=1, \dots, k$ and $R$, the radius of influence, we define
\[ I_i(R) = \left\{ \begin{array}{ll}
                 1 & \mbox{if $R=nk+i$ for some $n \in \bbN$ }\\
                 0 & \mbox{otherwise}.
                \end{array}
                \right. \]
Besides, we define
\[ \underline{x}_{n,i} := (\prod_{j=1}^{k} d_j)^{n}\prod_{j=1}^{i}d_{(j)} \hbox{ for } i \not= 0,\  \underline{x}_{n,0} := (\prod_{j=1}^{k} d_j)^{n} \hbox{ and } \underline{x}_{-1,i}:=0 \]
\noindent and
\[
\bar{x}_{n,i} := (\prod_{j=1}^{k} d_j)^{n}\prod_{j=1}^{i}d_{(k+1-j)} 
\hbox{ for } i \not= 0 \hbox{ and } \bar{x}_{n,0} := (\prod_{j=1}^{k} d_j)^{n}.
\]
and
\[h_i(R) =
\left[\sum_{m=0}^{\lfloor \frac{R-i}{k} \rfloor -1} \sum_{j=0}^{k-1} 
(\underline{x}_{m,j})^{-1} + \sum_{j=0}^{i-1} 
(\underline{x}_{\lfloor \frac{R-i}{k} \rfloor,j})^{-1}\right]G^R.
\]
\end{defn}

Analogously to definition \ref{D: RootedTree}, we consider the \cpp\ on $\bbT_{\tilde{d}}^+$. Relations analogous to~(\ref{E: NotEqual}) and~(\ref{E: Equal}) also holds between $\bbT_{\tilde{d}}$ and $\bbT_{\tilde{d}}^+$.

\begin{teo}
\label{T: CSPAH}
Consider the \cpp\ on $\bbT_{\tilde{d}}^+$ with radius of influence $R$
\begin{enumerate}
\item[\textit{(I)}] If  \begin{displaymath}
\sum_{i=0}^{k-1}c_i \mathbb{E}(G^R I_i(R)) > 1 + p_0
\end{displaymath}
 then, $\bbP_+(V) > 0,$
\item[\textit{(II)}] If
\begin{displaymath}
\sum_{i=0}^{k-1}\bar{c_i} \mathbb{E} \displaystyle
\left(h_i(R) I_i(R) \right) \leq 1
\end{displaymath}
then, $\bbP_+(V) = 0.$
\end{enumerate}
\end{teo}

\begin{cor}
\label{C: CSPAH}
Consider the \cpp\ on $\bbT_{d}^+$ (the $d-$dimensional rooted homogeneous tree) with radius of influence $R$
\begin{enumerate}
\item[\textit{(I)}] If $(1-p_0) d > 1$ then, $\bbP_+(V) > 0,$
\item[\textit{(II)}] If $(1-p_0) d \le 1$ and $\mathbb{E}(d^R) > 1 + p_0$ then, $\bbP_+(V) > 0,$
\item[\textit{(III)}] If $\mathbb{E}(d^R) \leq 2 - \frac{1}{d}$ then, $\bbP_+(V) = 0.$
\end{enumerate}
\end{cor}

Let $\rho$ and $\psi$ be, respectively, the smallest non-negative root of the equations
\begin{align}
\label{EQ: rho}
& \sum_{i=0}^{k-1}\mathbb{E}(\rho^{c_iG^R} I_i(R)) + (1-\rho) p_0 = \rho, \\
\label{EQ: psi}
& \sum_{i=0}^{k-1} \mathbb{E}(\psi^{\lfloor \bar{c_i}h(R)\rfloor} I_i(R)) = \psi,
\end{align}

\begin{teo}
\label{T: SobrevivenciaTd+}
Consider the \cpp\ on $\bbT_{\tilde{d}}^+.$
Then,
\begin{displaymath}
1 - \rho \leq \bbP_+(V) \leq 1 - \psi.
\end{displaymath}
\end{teo}

\begin{teo}
\label{T: ViveTd}
For the \cpp\ on $\bbT_{\tilde{d}}$ with radius of influence $R$, it holds that
\begin{equation*}
\label{E: ViveTd2}
1 - \sum_{i=0}^{k-1} \mathbb{E}\displaystyle
\left(\rho^{M_R({\mathcal O})} I_i(R) \right) \leq 
\bbP(V) \leq 1 - \sum_{i=0}^{k-1} \mathbb{E}\displaystyle
\left(\psi^{|T_R^{\mathcal O}|} I_i(R) \right).
\end{equation*}
\end{teo}

\begin{cor}
\label{C: ViveTd}
For the \cpp\ on $\bbT_d$ (the $d-$dimensional homogeneous tree) 
with radius of influence $R$, it holds that
\begin{equation*}
\label{E: ViveTdc}
1 - \displaystyle \left(1 -
\rho^{\frac{d+1}{d}}\right)p_0 -
\mathbb{E}\displaystyle \left(\rho^{\frac{(d+1)}{d}d^{R}}\right)
\leq \bbP(V) \leq 1 - \mathbb{E}\displaystyle
\left(\psi^{\frac{(d+1)}{d-1}(d^{R}-1)}\right)
\end{equation*}
where $\rho$ and $\psi$ are the smallest non-negative root of the equations~(\ref{EQ: rho}) and~(\ref{EQ: psi}).
\end{cor}

\begin{exa}
\label{E: G1-p}
Consider a \cpp\ in $\bbT_{\tilde{d}}, \tilde{d} = (4,9)$ assuming $R \sim {\mathcal G}(1-p).$ 
From Theorem~\ref{T: CSPAH} and equation~(\ref{E: Equal})
\[
0.078542 \leq \inf\{p: \bbP(V)>0\} \leq 0.097374.
\]
\end{exa}

\begin{exa}
\label{E: G2-p}
Consider a \cpp\ in $\bbT_{\tilde{d}}$, with $\tilde{d} = (12, 15, 16)$. Assuming
$R \sim {\mathcal B}(3, 0.1),$ from Theorem~\ref{T: ViveTd} we have,
\[
0.266557 \leq \bbP(V) \leq 0.266894.
\]

\end{exa}

\section{Spherically Symmetric Trees}
\label{S: Spherically Symmetric Trees}

\begin{defn}
We say that a tree, $\bbT_S$, is \textit{spherically symmetric}, if any pair of vertices at the same distance from the origin, have the same degree.
\end{defn}

Note that periodic trees are a subclass of spherically symmetric tree and therefore the results will also apply to periodic trees. In the previous section we obtained stronger results using the particular properties of periodic trees.

From definition~\ref{D: DefinicaoCPP} we consider the \cpp\ on $\bbT_S$.

\begin{defn}
Let us define for a tree ${\mathbb T}$
\begin{displaymath}
{\textrm{dim\ inf\ }  \partial {\mathbb T}} : = \lim_{n \rightarrow \infty}
\min_{v \in \mathcal{V}} \frac{1}{n} \ln M_n(v).
\end{displaymath}
\end{defn}

Observe that 
\[ \textrm{dim\ inf\ } \partial {\mathbb T}_d = \ln d.\]

\begin{teo}
\label{T: esferic}
For a \cpp\ in $\bbT_S$ and $R$, the radius of influence, $\bbP(V) > 0$ if
\begin{equation*}
\label{E: esferic}
\lim_{n \rightarrow \infty}\sqrt[n]{\rho_n} > 
e^{-\textrm{dim\ inf\ } \partial {\mathbb T}_S}
\end{equation*}
where
\begin{equation*}
\rho_n := \prod_{k=0}^{n-1} [1-\prod_{i=0}^{k}\bbP(R < i+1)].
\end{equation*}
\end{teo}

Lemma~\ref{L: aes2} shows that $\rho_n$ is as a lower bound of the probability that the process
starting from any vertex $v$ reaches the vertices at  $\partial T^v_n$,

\begin{cor}
\label{C: esferic1}
For a \cpp\ in $\bbT_S$ and $R$, a radius of influence satisfying $\bbP(R \leq k) = 1$ for some $k \in \bbN$, $\bbP(V) > 0$ if

\begin{equation*}
\label{E: esferic1}
\textrm{dim\ inf\ } \partial {\mathbb T}_S > \ln \displaystyle \left[\frac{1}{1 - \prod_{j=1}^{k}\bbP(R < j)} \right].
\end{equation*}
\end{cor}

\begin{cor}
\label{C: esferic2}
For a \cpp\ in $\bbT_S$ and $R$, a radius of influence satisfying
\begin{displaymath}
\bbP(R = k) = \frac{Z_{\alpha}}{(k+1)^{\alpha}}, \ k=1,2,\dots
\end{displaymath}
if dim\ inf\  $\partial {\mathbb T}_S > 0,$ then $\bbP(V) > 0.$
\end{cor}

\begin{exa}
\label{E: B(p)}
Consider  a \cpp\ in $\bbT_S$ with $R \sim {\mathcal B}(p).$
\begin{itemize}
\item If dim\ inf\ $\partial {\mathbb T}_S > -\ln p$ then, $\bbP(V)>0,$
\item If $\bbT_S = \bbT_{\tilde{d}}$ and $ G(\tilde{d}) > \frac{1}{p}$ then, $\bbP(V)>0.$
\end{itemize}

\end{exa}

\section{Galton-Watson Branching Trees}
\label{S: Galton-Watson Branching Trees}

\subsection{Non Homogeneous Galton-Watson Branching Trees}
\label{S: Non Homogeneous G-W-B Trees}

Consider a supercritical Galton-Watson branching process starting from a single progenitor 
such that each individual whose distance from the progenitor is $n$ has a random number of 
offspring (independet of everything else) with generating function $f_n(s) = \sum_{k=0}^{\infty} q_n(k)s^k$. 

Let us define $F=\{(f_n, d_n)\}_{n \in \bbN}$ where $d_n = f^{\prime}_n(1) \in (0, \infty) $.  This Galton-Watson branching process yields a random family tree $\bbT_{F}$. We are 
particularly interested in a supercritical Galton-Watson tree, on the event of non extinction (infinite trees). 
A sufficient condition for that is $\liminf_{n \to \infty} d_n > 1 $.

\begin{defn}
For a supercritical Galton-Watson tree on $\bbT_F$, let us define
\begin{displaymath}
D(\bbT_F) : = \lim_{n \rightarrow \infty}
\min_{i \in \mathbf{N}} \frac{1}{n} \ln \left[\prod_{l=i}^{i+n-1}d_l \right].
\end{displaymath}
\end{defn}

In particular, if $F$ is a sequence of generating functions of degenerated random 
variables $\{X_n\}_{n \ge 0}$ such that $X_n=a_n$ we have that $\bbT_F$ equals
to a spherically symmetric tree $\bbT_S$ with probability 1. 
Then, with probability 1
\begin{displaymath}
D(\bbT_F) = {\textrm{dim\ inf\ }  \partial {\mathbb T_S}}.
\end{displaymath}

\begin{teo}
\label{T: NHGWT}
For a \cpp\ on $\bbT_F$ with a radius of influence $R$, $\bbP(V) > 0$ if
\begin{equation*}
\label{E: esferic}
\lim_{n \rightarrow \infty}\sqrt[n]{\rho_n} >
e^{-D(\bbT_F)}
\end{equation*}
where
\begin{equation*}
\rho_n := \prod_{k=0}^{n-1} [1-\prod_{i=0}^{k}\bbP(R < i+1)].
\end{equation*}
\end{teo}

\subsection{Homogeneous Galton-Watson Branching Trees}
\label{S: Homogeneous G-W-B Trees}

Consider a supercritical Galton-Watson branching process starting from a single progenitor 
such that each individual has a random number of offspring (independet of everything else) 
whose average is $d > 1$. This process yields a random infinite family tree, known as a 
supercritical Galton-Watson tree $\bbT_{F}$, where $d_n = d$ for all ${n \in \bbN}$, on the event of non extinction.

\begin{teo}
\label{T: HGWT}
Consider the \cpp\ on a homogeneous supercritical Galton-Watson branching tree with radius of influence $R$.
\begin{enumerate}
\item[\textit{(I)}] If $(1-p_0) d > 1$ then, $\bbP[V] > 0,$
\item[\textit{(II)}] If $(1-p_0) d \le 1$ and $\mathbb{E}(d^R) > 1 + p_0$ then, $\bbP[V] > 0,$
\item[\textit{(III)}] If $\mathbb{E}(d^R) \leq 2 - \frac{1}{d}$ then, $\bbP[V] = 0.$
\item[\textit{(IV)}] For $\mathbb{E}(d^R) < 2 - \frac{1}{d}$ we have
\begin{displaymath}
 \frac{1 }{1-\mathbb{E}\displaystyle \left (d^R \right) + p_0} \leq \mathbb{E}(|I|)  \leq \frac{d-1}{2d - 1 - d\mathbb{E}\displaystyle \left (d^R \right)}.
\end{displaymath}
\end{enumerate}
\end{teo}

\section{Proofs}
\label{S: Proofs}

\subsection{Homogeneous Trees}
\hfill

\noindent
\textit{Proof of Theorem~\ref{T: THT1}}\\
Let us define now two auxiliary branching process. For the first,
$\{\mathcal{X}_n\}_{n \in \bbN}$, each individual has a number of
offspring distributed as the random variable $X$, assuming values
in $\{0, d, d^2, \dots \}$ such that
\[ \bbP[X=0]= p_o, \bbP[X=d]= p_1, \cdots, \bbP[X=d^k]= p_k
\text{ for all } k=1, 2, \dots \]

In the second auxiliary process, $\{\mathcal{Y}_n\}_{n \in \bbN}$,
each individual has a number of offsprings distributed as the random
variable $Y$, assuming values in $\{ 0, d, d+d^2, \dots, 
\sum_{i=1}^k d^i, \dots \}$ such that
\[ \bbP[X=0]= p_o, \bbP[X=d]= p_1, \cdots, \bbP[X=\sum_{i=1}^k d^i]= p_k
\text{ for all } k=1, 2, \dots \]

These two processes provide convenient lower bounds ($\{\mathcal{X}_n\}_{n \in \bbN}$) and upper bounds ($\{\mathcal{Y}_n\}_{n \in \bbN}$) for our process. Suppose that 
$R_v=r$ for a fixed site $v$. Then the set of vertices activated by $v$ is $T_r^v$, whose cardinality is $\sum_{i=1}^r d^i$ vertices.
The activation process will go on. The process $\{\mathcal{X}_n\}_{n \in \bbN}$ will only count on those $d^k$ which are at distance $r$
from $v$ (the set $\partial T_r^v$). By the other side, the process $\{\mathcal{Y}_n\}_{n \in \bbN}$ counts activation that will be made by all of them ($T_r^v$), in addition to disregarding the fact that some vertice will experience multiple activations from sites belonging to $T_r^v$.

For these processes the average number of offsprings are respectively
$\mu_X = \mathbb{E}\displaystyle \left(d^R \right) - p_0 $ and $\mu_Y = \frac{d}{d-1}\displaystyle 
\left[\mathbb{E}\displaystyle \left(d^R \right) - 1 \right] $.  As $\mu_X < 1$ and $\mu_Y < 1$ by
hypothesis, the expected values for the total number of individuals are respectively
\[ \frac{1}{1-\mu_X}  = \frac{1}{1+p_0-\mathbb{E}\displaystyle \left(d^R \right)} \] and
\[ \frac{1}{1-\mu_Y} = \frac{d-1}{2d-1-d\mathbb{E}\displaystyle \left(d^R \right)}.\]

Using the fact that the root has degree $d+1$ we can modify the processes $\{\mathcal{X}_n\}_{n \in \bbN}$ and $\{\mathcal{Y}_n\}_{n \in \bbN}$ such that the offspring distributions for the first generation are res\-pectively
\begin{align*}
&\mathbb{P}[X=0] = p_0, \\
&\mathbb{P}[X=(d+1)d^{k-1}] = p_k \hbox{ for } k =1,2, \dots
\end{align*}

and
\begin{align*}
&\mathbb{P}\Big[Y = \frac{(d+1)(d^k-1)}{d-1}\Big] = p_k \hbox{ for } k = 0,1,2, \dots
\end{align*}

For these modified processes the total expected number of individuals are respectively

\begin{displaymath}
\mathbb{E}(|I_x|) = \sum_{k=1}^{\infty}\displaystyle \left(\frac{(d+1)d^{k-1}}{1 +p_0 - \mathbb{E}\displaystyle \left(d^R \right) } +
1 \right)p_k + p_0 =  \frac{d + \mathbb{E}\displaystyle \left (d^R \right) - p_0}{d(1-\mathbb{E}\displaystyle \left (d^R \right) + p_0)}
\end{displaymath}
and
\begin{align*}
\mathbb{E}(|I_y|) &= \sum_{k=0}^{\infty}\left(  \displaystyle \left[\frac{(d+1)(d^k-1)}{(d-1)} \right]\displaystyle \left[\frac{(d-1)}{2d-1-d\mathbb{E}\displaystyle \left (d^R \right)}\right]+1 \right) p_k \\
&= \frac{\mathbb{E}\displaystyle \left (d^R \right)+d-2}{2d - 1 - d\mathbb{E}\displaystyle \left (d^R \right)}.
\end{align*}
Since the reasonings that justified the lower and upper bounds at the beginning of the proof are valid with this modification, we have that
$\mathbb{E}(|I_x|) \leq \mathbb{E}(|I|) \leq \mathbb{E}(|I_y|)$ and the result follows. \qed

\subsection{Periodic Trees}
\hfill

Consider a  $k-$periodic tree whose degrees are
$d_{1}+1, d_{2}+1, \cdots, d_{k}+1$ and for $ i= 1, \dots, k-1 $

\[ J_i = \{(j_1, \dots, j_k), 1 \leq j_1 < j_2 < \cdots < j_i \leq k \}. \]

Let us define for $n \in \bbN$
\begin{align*}
A_{nk} &= \{(\prod_{j=1}^{k} d_j)^{n} \}, \\ 
A_{nk+i} &= \{( \prod_{j=1}^{k} d_j)^{n}\prod_{l=1}^{i}d_{j_l}, (j_1, \dots, j_i) \in J_i \}  \hbox{ for } i= 1, \dots, k-1 .
\end{align*}

We claim that for all $n \in \bbN$, $k \in \bbN$ and $v \neq \mathcal{O}$ that
\begin{align}
\label{Claim 1}
\min A_{nk+i} &= \underline{x}_{n,i}, \\
\label{Claim 2}
\max A_{nk+i} &= \bar{x}_{n,i}, \\
\label{Claim 3}
M_{nk+i} (v)&\in A_{nk+i}.
\end{align}

Let
\begin{align*}
y_{n,i} := \sum_{m=0}^{n-1} \sum_{j=0}^{k-1} 
(\underline{x}_{m,j})^{-1} + \sum_{j=0}^{i-1} 
(\underline{x}_{m,j})^{-1}
\end{align*}

\begin{lem} 
\label{L: ppakpii}
Consider a $k-$periodic tree whose degrees are $d_1+1, d_2+1, \cdots,
d_k + 1$, $d_i\geq$ 2 for all $i = 1,2, \cdots k$. Consider a vertex $v \neq \mathcal{O}$.
Then
\begin{displaymath}
|T_{nk+i}^v|  \leq \lfloor y_{n,i} \cdot  \bar{x}_{n,i} \rfloor.
\end{displaymath}
\end{lem}

\noindent
\textit{Proof of Lemma~\ref{L: ppakpii}}\\
Consider first the following set up: $R=k, d({\mathcal O},v) = mk$ for some $m \in \bbN$ and $\bbT_{\tilde{d}}$ such that $d_i = d_{(i)}$
for all $i = 1, \cdots, k$. Then
\begin{displaymath}
 \displaystyle {
|T_{nk+i}^v| = \left|\bar{x}_{n,i} 
+ \frac{\bar{x}_{n,i}}{d_{k}} + \frac{\bar{x}_{n,i}}{d_{k}d_{k-1  }}
+ \cdots
+ \frac{\bar{x}_{n,i}}{\prod_{j=1}^{k}d_{j}}\right|}.
\end{displaymath}

Consider now the case where $R=nk$ and $d({\mathcal O},v) = mk$ for $n, m \in \bbN$ and $\bbT_{\tilde{d}}$ such that $d_i = d_{(i)}$
for all $i = 1, \cdots, k$. Then
\begin{align*}
|T_{nk+i}^v| &=\left |  \bar{x}_{n,i} 
+ \frac{\bar{x}_{n,i}}{d_{k}} + \frac{\bar{x}_{n,i}}{d_{k}d_{k-1}}
+ .
+ \frac{\bar{x}_{n,i}}{\prod_{j=1}^{k}d_{j}} \right | \\
&+ \left |
\frac{\bar{x}_{n,i}}{(\prod_{j=1}^{k}d_{j})d_{k}} 
+ \cdots
+ \frac{\bar{x}_{n,i}}{(\prod_{j=1}^{k}d_{j})^2} \right | + \cdots \\
&+ \left |
\frac{\bar{x}_{n,i}}{(\prod_{j=1}^{k}d_{j})^{n-1}d_{k}} 
+ \cdots + \frac{\bar{x}_{n,i}}{(\prod_{j=1}^{k}d_{j})^n}
\right |.
\end{align*}

Observe now that on any tree, for any $v_r$ such that $d(\mathcal{O}, v_r)= r$
\[ |T_{m}^{v_r}| = \sum_{j=1}^{m} M_j(v_r) =  M_m(v_r) + \sum_{j=1}^{m-1} M_m(v_r). \left [ \prod_{i=1}^{j} M_1(v_{r+m-i}) \right ]^{-1}
\]

Now consider only $R=nk+i$ for $n \in \bbN$ and $i = 1, \cdots, k-1$. So, from~(\ref{Claim 1}), (\ref{Claim 2}) and~(\ref{Claim 3}), it follows that
\begin{align*}
|T_{nk+i}^v| &\leq \Big\lfloor  \bar{x}_{n,i} 
+ \frac{\bar{x}_{n,i}}{d_{(1)}} + \frac{\bar{x}_{n,i}}{d_{(1)}d_{(2)}}
+ \cdots 
+ \frac{\bar{x}_{n,i}}{\prod_{j=1}^{k}d_{(j)}} + \cdots \\
&+ \frac{\bar{x}_{n,i}}{(\prod_{j=1}^{k-1}d_{(j)})^{n}d_{(1)}}
+ \cdots
+ \frac{\bar{x}_{n,i}}{(\prod_{j=1}^{k}d_{(j)})^n{\prod_{j=1}^{i}d_{(j)}}} \Big\rfloor \\
&= \displaystyle \Big\lfloor
y_{n,i} \cdot \bar{x}_{n,i} \Big\rfloor.
\end{align*} \qed

Let us define two auxiliary branching process, being the first one
$\{\mathcal{X}_n\}_{n \in \bbN}$. This process is defined by a random variable 
$ X,$ assuming values in 
$\{ \underline{x}_{n,i}, i=0, \dots, k-1, \hbox{ and } n=0, 1, \dots, (n,i) \not = (0,0)\} \cup \{0\}$
such that
\begin{align*}
\mathbb{P}[X=0] &=: p_0, \\
\mathbb{P}[X=\underline{x}_{n,i}] &=: p_{nk+i} \hbox{ for } i=0, \dots, k-1, \hbox{ and } n=0, 1, \dots, (n,i) \not = (0,0)
\end{align*}

Its expected value is given by the following lemma
\begin{lem} 
\label{L: esp1}
\begin{align*}
\label{E: Minorante}
\mathbb{E}[X] = \sum_{i=0}^{k-1}c_i\mathbb{E}\displaystyle \left[G^{R}I_{i}(R)
\right] - p_0
\end{align*}
\end{lem}
\noindent
\textit{Proof of Lemma~\ref{L: esp1}}\\
\begin{align*}
\mathbb{E}(X) &= \sum_{n=1}^{\infty}\underline{x}_{n,0}p_{nk}+
\sum_{i=1}^{k-1} \sum_{n=0}^{\infty}\underline{x}_{n,i}p_{nk+i} \\
&= \sum_{i=0}^{k-1}c_i
\sum_{n=0}^{\infty}\prod_{j=1}^{k}(\sqrt[k]{d_j})^{nk+i}p_{nk+i}-p_0 \\&= \mathbb{E}[G^{R}I_0(R)] +
\sum_{i=1}^{k-1}c_i\mathbb{E}[G^{R}I_i(R)] - p_0
\\ &= \sum_{i=0}^{k-1}c_i\mathbb{E}[G^{R}I_i(R)] -
p_0.
\end{align*} \qed

\noindent
and its probability generating function is given by

\begin{lem} 
\label{L: fgp1}
\begin{equation*}
\label{G: fgp}
\varphi_{X}(s) = \sum_{i=0}^{k-1}\mathbb{E}\displaystyle \left[s^{c_i G^{R}}I_{i}(R)
\right]  + (1-s)p_0.
\end{equation*}
\end{lem}

\noindent
\textit{Proof of Lemma~\ref{L: fgp1}}
\begin{align*}
\varphi_{X}(s) & = p_{0} +
\sum_{n=1}^{\infty}s^{\underline{x}_{n,0}}p_{nk} + \sum_{i=1}^{k-1}\sum_{n=0}^{\infty}s^{\underline{x}_{n,i}}p_{nk+i}
\\& = p_{0} + \sum_{n=1}^{\infty}s^{G^{nk}}p_{nk} + \sum_{i=1}^{k-1}\sum_{n=0}^{\infty}s^{c_iG^{nk+i}}p_{nk+i}
\\&  = 
p_{0} - s p_{0} + \sum_{i=0}^{k-1}\sum_{n=0}^{\infty}s^{c_iG^{nk+i}}p_{nk+i} 
\\&= (1- s)p_{0} +
\sum_{i=0}^{k-1}\mathbb{E}\displaystyle
\left[s^{c_i G^R}I_i (R) \right].
\end{align*}

The second auxiliary process is $\{\mathcal{Y}_n\}_{n \in \bbN}$, a branching process defined by a random variable $ Y,$ assuming values on 
$\{ \lfloor y_{n,i} \bar{x}_{n,i} \rfloor, i=0, \dots, k-1, \hbox{ and } n=0, 1, \dots\} $
such that
\begin{align*}
&\mathbb{P}\Big[Y = \lfloor y_{n,i} \bar{x}_{n,i}
 \rfloor \Big] = p_{nk + i} \hbox{ for } i = 0,1, \dots, k-1 \hbox{ and } n=0,1,\dots
\end{align*}
Its expected value satisfies
\begin{lem} 
\label{L: esp2}
\begin{equation*}
\label{E: Majorante}
\mathbb{E}[Y] \leq  \sum_{i=0}^{k-1}\bar{c_i}\mathbb{E}\displaystyle \left[ \displaystyle h_i(R) I_i(R) \right].
\end{equation*}
\end{lem}
\noindent
\textit{Proof of Lemma~\ref{L: esp2}}
\begin{align*}
\mathbb{E}(Y)  \leq & 
\sum_{i=0}^{k-1} \sum_{n=0}^{\infty} y_{n,i}
 \bar{x}_{n,i} p_{nk + i} \\
=& \sum_{i=0}^{k-1} \sum_{n=0}^{\infty}
\left[\sum_{m=0}^{n-1} \sum_{j=0}^{k-1} 
(\underline{x}_{m,j})^{-1} + \sum_{j=0}^{i-1} 
(\underline{x}_{m,j})^{-1}\right]
 \bar{x}_{n,i} p_{nk + i} \\
=&\sum_{i=0}^{k-1}
\bar{c}_i\sum_{n=0}^{\infty}
\left[\sum_{m=0}^{n-1} \sum_{j=0}^{k-1} 
(\underline{x}_{m,j})^{-1} + \sum_{j=0}^{i-1} 
(\underline{x}_{m,j})^{-1}\right]
\prod_{j=1}^{k}(\sqrt[k]{d_j})^{nk+i} p_{nk + i} \\
=& \sum_{i=0}^{k-1}\bar{c_i}
\mathbb{E}\displaystyle \left[ \displaystyle h_i(R) I_i(R) \right].
\end{align*} \qed

\noindent
and its probability generating function is given by
\begin{lem} \label{L: fgp2}
\begin{align}
\label{G: fgp2}
\varphi_{Y}(s) =  \sum_{i=0}^{k-1}\mathbb{E}\displaystyle \left[s^{
\displaystyle \left\lfloor \bar{c}_i h_i(R)\right \rfloor }I_{i}(R)
\right].
\end{align}
\end{lem}
\noindent
\textit{Proof of Lemma~\ref{L: fgp2}}
\begin{align*}
\varphi_{Y}(s) =& \sum_{i=0}^{k-1} \sum_{n=0}^{\infty}s^{\lfloor
y_{n,i}\bar{x}_{n,i}
 \rfloor}p_{nk + i}
\\=& \sum_{i=0}^{k-1} \sum_{n=0}^{\infty}s^{\lfloor
y_{n,i}G^{nk+i}\bar{c}_{i}
 \rfloor}p_{nk + i} \\=& \sum_{i=0}^{k-1}\mathbb{E}\displaystyle \left[s^{
\displaystyle \left\lfloor \bar{c}_i h_i(R)\right \rfloor }I_{i}(R)
\right].
\end{align*} \qed

\noindent
\textit{Proof of Theorem~\ref{T: CSPAH}}\\
By a coupling argument one can see that our process dominates (by~(\ref{Claim 1}) and~(\ref{Claim 3}))
$\{\mathcal{X}_n\}_{n \in \bbN}$. This process survives as long as
$\mathbb{E}[X]> 1$. Therefore from Lemma~\ref{L: esp1} our process survives
if
\[ \sum_{i=0}^{k-1}c_i\mathbb{E}\displaystyle \left[G^{R}I_{i}(R)
\right] > 1 + \mathbb{P}(R = 0),\] 
\noindent
proving \textit{(I)}.

By the other side, also by a coupling argument, our process is do\-mi\-na\-ted (by~(\ref{Claim 2}) and~(\ref{Claim 3})) by
$\{\mathcal{Y}_n\}_{n \in \bbN}.$ That process dies out provided $\mathbb{E}[Y] \leq 1$
therefore from Lemma~\ref{L: esp2} our process dies out if
\[ \sum_{i=0}^{k-1}\bar{c_i} \mathbb{E} \displaystyle
\left(h_i(R) I_i(R) \right) \leq 1,\]
\noindent 
proving \textit{(II)}.\qed

\noindent
\textit{Proof of Theorem~\ref{T: SobrevivenciaTd+}}

In order to find the extinction probability of $\{\mathcal{X}_n\}_{n \in \bbN}$
(Grimmett and Stirzaker(~\cite[p.173]{GrimmettStirzaker}),
let us consider the smallest non-negative root of the equation $\rho = \varphi_{X}(\rho).$
Therefore from Lemma~\ref{L: fgp1}
\[ \sum_{i=0}^{k-1}\mathbb{E}\displaystyle
\left[\rho^{c_{i}G^{R}}I_{i}(R) \right] + (1 - \rho)p_0
= \rho \]
and by construction of the processes, as $ \bbP_+[V^c] \leq \rho, $ we have that
\[ 1- \rho \leq \bbP_+(V) .\]
In order to find the extinction probability of $\{\mathcal{Y}_n\}_{n \in \bbN}$
(Grimmett and Stirzaker~\cite[p.173]{GrimmettStirzaker}),
let us consider the smallest non-negative root of the equation $\psi = \varphi_{Y}(\psi).$
Therefore from Lemma~\ref{L: fgp2}
\[ \sum_{i=0}^{k-1}\mathbb{E}\displaystyle \left[\psi^{
\displaystyle \left\lfloor \bar{c}_i h_i(R)\right \rfloor }I_{i}(R)
\right]) = \psi \]
and by the construction of the processes, as $ \bbP_+[V^c] \geq \psi$, we have that
\[ \bbP_+(V) \leq 1 - \psi.\]
\qed

\noindent
\textit{Proof of Theorem~\ref{T: ViveTd}}

Observe that except for the root, all vertices see towards
infinity a tree like $\bbT_{\tilde{d}}^+.$ So, assuming $R_{\mathcal O}=nk+i$ the probability
for the process to survive is greater or equal than the probability of the
process to survive from at least one of the $M_{nk+i}{(\mathcal O)}$ trees that have
as root the furthest infected vertices. Now note that, still assuming
$R_{\mathcal O}=nk+i$, the probability for the process to survive on $\bbT_{\tilde{d}}$  is smaller
or equal than the probability for the process to survive from at least one
of the $|T_{nk+i}^{\mathcal O}|$ vertices which are in the radius of
influence ($R_{\mathcal O}$) of the origin of the tree as if each one had its own tree.
Then
\[
 \bbP(V | R_{\mathcal O} = nk+i)  \geq 1 - (1 - \bbP_+(V))^{M_{nk+i}{(\mathcal O)}} \geq 1 - \rho^{M_{nk+i}{(\mathcal O)}}
\]
and
\[
 \bbP(V | R_{\mathcal O} = nk+i)  \leq 1 - (1 - \bbP_+(V))^{|T_{nk+i}^{\mathcal O}|} \leq 1 - \psi^{|T_{nk+i}^{\mathcal O}|}.
\]
Then,
\begin{align*}
\mathbb{P}(V) =&  
\sum_{i=0}^{k-1}\sum_{n=0}^{\infty}\mathbb{P}(V | R_{\mathcal O} =
nk+i)p_{nk + i} 
\\ \geq & 
\sum_{i=0}^{k-1}\sum_{n=0}^{\infty}[1 -
\rho^{M_{nk+i}{(\mathcal O)}}]p_{nk + i}
\\
= & 
1-\sum_{i=0}^{k-1}\mathbb{E}[\rho^{M_{R}{(\mathcal O)}}I_i(R)]
\end{align*}
and
\begin{align*}
\mathbb{P}(V) & = \sum_{i=0}^{k-1}\sum_{n=0}^{\infty}\mathbb{P}(V | R_{\mathcal O} =
nk+i)\mathbb{P}(R_{\mathcal O} = nk+i)  \\
& \leq \sum_{i=0}^{k-1}\sum_{n=0}^{\infty}\displaystyle \left[1 -
\psi^{|T_{nk+i}^{\mathcal O}|} \right]p_{nk + i}
\\ &= 1 - \sum_{i=0}^{k-1}\mathbb{E}\displaystyle
\left(\psi^{|T_{R}^{\mathcal O}|} I_i(R) \right).
\end{align*}
\qed

\noindent
\subsection{Spherically Symmetric Trees}
\hfill

Suppose we have a set of independent random variables
$\{R_v\}_{\{ v \in {\mathcal V}(\bbT_S) \}}$ distributed as $R$. 
Assume $\bbP(R=0)<1.$

For $u \leq v \in {\mathcal V}(\bbT_S)$, consider the event
\[ V_{u,v}: \textit{Process starting from $u$ reaches $v$}.\]

For a fixed integer $n,$ let $ X_0^n = \{{\mathcal O}\}.$ Besides, for
$ j=1,2, \dots $ consider
\[ X_j^n = \bigcup_{u \in X_{j-1}^n} \{ v \in \partial
T_n^u : V_{u,v} \textit{ occurs } \}.\]
Again, for all $ j=1,2, \dots $ consider
\[Z_j^n = | X_j^n |.\]

So, for all fixed positive integer $n$, $\{Z_j^n\}_{j \geq 0}$ is a branching process dominated
by the number of vertices $v \in \partial T^{\mathcal O}_{jn}$ which are activated.

\begin{lem}
\label{L: aes1}
Consider $n$ fixed. For $\mu_j,$ the mean number of offspring
of one individual of generation $j$ for the process
$\{Z_j^n\}_{j \geq 0}$, it holds that
\begin{displaymath}
\mu_j := \mu_j^n  = M_n(u) \rho_{j}^{(n)},
\end{displaymath}
where $\rho_{j}^{(n)} = \bbP (V_{u,v})$, for any fixed pair $u \leq v$ such that $d({\mathcal O},u) =jn$
and $d({\mathcal O},v)=(j+1)n.$
\end{lem}
\noindent 
\textit{Proof of Lemma~\ref{L: aes1}}\\
For fixed $j$ and $n$, consider for some $u$ such that $d({\mathcal O},u) =jn$,
$\partial T_n^u = \{ v_1, v_2, \dots , v_{M_n(u)} \}$. So we can write the
number of offspring of $u$ as $\sum_{i=1}^{M_n(u)}I_{\{V_{u,v_i} \}}.$ Taking expectation finishes the proof. \qed

\begin{lem} 
\label{L: aes2}
Consider $n$ fixed and $\rho_{j}^{(n)} = \bbP (V_{u,v})$, for any fixed pair $u \leq v$ such that
$d({\mathcal O},u) =jn$ and $d({\mathcal O},v)=(j+1)n,$
\begin{displaymath}
\rho_{j}^{(n)} \geq
\prod_{k=0}^{n-1}[1-\prod_{i=0}^{k}\bbP (R < i+1)].
\end{displaymath}
\end{lem}
\noindent 
\textit{Proof of Lemma~\ref{L: aes2}}\\
For any fixed pair $u \leq v$ such that $d({\mathcal O},u) =jn$ and $d({\mathcal O},v)
=(j+1)n$ we have that

\begin{displaymath}
V_{u,v} = \bigcap_{k=0}^{n-1}\displaystyle
\left[\bigcup_{i=0}^{k}\{ R_{u(i)} \geq k+1-i \}\right]
\end{displaymath}
\noindent
where $u(i)$ is the vertex from the path conecting $u$ to $v$
such that $d({\mathcal O},u(i))=jn+i$. From this follows

\begin{align*}
\rho_{j}^{(n)} = & \ \bbP \displaystyle \left (
\bigcap_{k=0}^{n-1}\displaystyle \left[\bigcup_{i=0}^{k}\{
R_{u(i)} \geq k+1-i \}\right] \right) \\ \geq & \
\prod_{k=0}^{n-1}\bbP \displaystyle \left(\bigcup_{i=0}^{k}\{
R_{u(i)} \geq k+1-i \}\right).
\end{align*}
The inequality is a consequence of the FKG inequality (N.Alon and J.Spencer~\cite[p.89]{AS}). \qed

\noindent
\textit{Proof of Theorem~\ref{T: esferic}}\\
Assume that dim inf $\partial {\mathbb T}_S > 0$. Then, for all $\alpha \in
(0,$dim inf $\partial {\mathbb T}_S$ ) there exists $N = N(\alpha)$ such that for all $n \geq N$
\begin{displaymath}
\min_{v \in \mathcal{V}}\frac{1}{n}\ln M_n(v) > \alpha
\end{displaymath}
where
\begin{displaymath}
M_n(v) \geq e^{\alpha n} \textrm { for all } v \in \mathcal{V}
\textrm { and } n \geq N.
\end{displaymath}

From Souza $\&$ Biggins ~(\cite[p.40]{SouzaBiggins}) a branching process in varying environments is
\textit{uniformly supercritical} if there exists constants $a > 0$ and $c > 1$
such that
\begin{equation*} 
\label{UniformlySupercritical}
 \prod_{k=i}^{j+i-1}\mu_k \geq ac^j, \textrm { for all } i \geq 0 \textrm { and } j \geq 0.
\end{equation*}
Observe that that condition holds if
\begin{displaymath}
\liminf_{j \rightarrow \infty} \mu_j > 1
\end{displaymath}
From Lemma~\ref{L: aes1} we have that for $n \geq N$
\begin{displaymath}
\liminf_{j \rightarrow \infty} \mu_j \geq e^{\alpha n}\rho_n =
(e^{\alpha}\sqrt[n]{\rho_n})^{n}
\end{displaymath}

Now note that we can write
\begin{displaymath}
Z_{j+1} = \sum_{i=1}^{Z_j}Y_{j,i}^n,
\end{displaymath}
where $Y_{j,i}^n$ are i.i.d. copies of $Y_j^n$, being the number of offspring from the $i-th$ individual of the $j-th$ generation. By considering
Lemma~\ref{L: aes1} we have for all $j$ that
\begin{displaymath}
\frac{Y_j^n}{\mu_j} \leq \frac{M_n(u)}{\mu_j} = \frac{1}{\rho_j^{(n)}} \leq (\bbP [R > 0])^{-n}
\end{displaymath}
where $\rho_{j}^{(n)} = \bbP (V_{u,v})$, for any fixed pair $u \leq v$ such that $d({\mathcal O},u) =jn$
and $d({\mathcal O},v)=(j+1)n.$

So, from Theorem 1 in Souza $\&$ Biggins~(\cite[p.40]{SouzaBiggins}), we conclude that the cone percolation process has a giant component with positive
probability if
\begin{displaymath}
\lim_{n \rightarrow \infty}
\sqrt[n]{\rho_n} > e^{-\alpha}.
\end{displaymath}
As this hold for every $\alpha \in (0,$dim inf $\partial {\mathbb T}_S$ ), the condition
\begin{displaymath}
\lim_{n \rightarrow \infty}\sqrt[n]{\rho_n} > e^{\textrm { - dim
inf }\partial {\mathbb T}_S}
\end{displaymath}
guarantees the survival of the process with positive probability. \qed

\noindent 
\textit{Proof of Corollary~\ref{C: esferic1}}\\
\begin{align*}
&\sqrt[n]{\prod_{i=0}^{n-1}[1-\prod_{j=0}^{i}\mathbb{P}(R <
j+1)]} = \\
& = [1 - \prod_{j=1}^{k}\mathbb{P}(R <
j)]\sqrt[n]{\frac{\prod_{i=0}^{k-1}[1-\prod_{j=0}^{i}\mathbb{P}(R
< j+1)]}{(1 - \prod_{j=1}^{k}\mathbb{P}(R < j))^k}} \\ 
&\rightarrow 1 - \prod_{j=1}^{k}\mathbb{P}(R < j), \textrm { when
} n \rightarrow \infty.
\end{align*}

\noindent
\textit{Proof of Corollary~\ref{C: esferic2}}\\
Observe that
\begin{align*}
\rho_n &\geq \mathbb{P}(R \geq n) =
\sum_{k=n}^{\infty}\frac{Z_{\alpha}}{(k+1)^{\alpha}}  \\
&\geq
\int_{n+1}^{\infty}\frac{Z_{\alpha}}{x^{\alpha}}dx =
\frac{Z_{\alpha}}{(\alpha - 1)(n+1)^{\alpha - 1}}
\end{align*}
The above inequalitty follows from the integral test.\\
Now observe that if $\textrm {dim inf }\partial {\mathbb T}_S > 0$, we have that
\begin{displaymath}
\lim_{n \rightarrow \infty} \sqrt[n]{\rho_n} \geq \lim_{n
\rightarrow \infty}\sqrt[n]{\frac{Z_{\alpha}}{(\alpha -
1)}\frac{1}{(n+1)^{\alpha - 1}}} = 1 > e^{- \textrm {dim inf
}\partial {\mathbb T}_S}
\end{displaymath}
Theorem~\ref{T: esferic} guarantess the desired result. \qed

\subsection{Galton-Watson Branching Trees}
\hfill

\subsubsection{Non Homogeneous Galton-Watson Branching Trees}
\hfill

\noindent
\textit{Proof of Theorem~\ref{T: NHGWT}}\\

Suppose we have a set of independent random variables
$\{R_{n,m}\}_{\{n, m \in \bbN \}}$ distributed as $R$.
Assume $\bbP(R=0)<1.$ For each tree $\bbT_f$ on $\bbT_F$ we associate each of
its existing vertices to a pair $u =  (n, m)$ so that $R_{n,m}$ is its radius of influence.
With this aim, $n$ stands for the distance from a set of $k(n)$ vertices to the tree progenitor
while $m=1, \cdots, k(n)$ stands for an enumeration on the set of the existing vertices at level $n$.

For each tree $\bbT_f$ on $\bbT_F$ and $u \leq v \in {\mathcal V}(\bbT_f)$, consider the event
\[ V_{u,v}: \textit{Process starting from $u$ reaches $v$}.\]

Let 
\[ \Omega = \{ (\bbT_f; \{r_{n,m}\}_{\{n, m \in \bbN \}} ); \bbT_f \in \bbT_F; \{r_{n,m}\}_{\{n, m \in \bbN \}} \in \bbN^{\bbN \times \bbN} \} \] 

Take $\omega = (\bbT_f; \{r_{n,m}\}_{\{n, m \in \bbN \}})$.
For a fixed integer $n,$ let $ X_0^n(\omega) = \{{\mathcal O}\}.$ Besides, for
$ j=1,2, \dots $ consider
\[ X_j^n(\omega) = \bigcup_{u \in X_{j-1}^n(\omega)} \{ v \in \partial
T_n^u(\omega) : I_{V_{u,v}}(\omega)=1 \}.\]
The definition for $\partial
T_n^u(\omega)$ is analogous to~(\ref{E: defPartialT}). 
Again, for all $ j=1,2, \dots $ consider
\[Z_j^n = | X_j^n |.\]

So, for all fixed positive integer $n$, $\{Z_j^n\}_{j \geq 0}$ is a branching process dominated
by the number of vertices $v \in \partial T^{\mathcal O}_{jn}$ which are activated.

\begin{lem}
\label{L: aes2}
Consider $n$ fixed. For $\mu_j,$ the mean number of offspring
of one individual of generation $j$ for the process
$\{Z_j^n\}_{j \geq 0}$, it holds that
\begin{displaymath}
\mu_j := \mu_j^n  = \left[\prod_{i=jn+1}^{jn+n}d_i \right] \rho_{j}^{(n)},
\end{displaymath}
where $\rho_{j}^{(n)} = \bbP (V_{u,v})$, for any fixed pair $u \leq v$ such that $d({\mathcal O},u) =jn$
and $d({\mathcal O},v)=(j+1)n.$
\end{lem}
\noindent
\textit{Proof of Lemma~\ref{L: aes2}}\\
For fixed $j$ and $n$, consider for some $u$ such that $d({\mathcal O},u) =jn$,
$\partial T_n^u = \{ v_1, v_2, \dots , v_{M_n(u)} \}$. So we can write the
number of offspring of $u$ as $\sum_{i=1}^{M_n(u)}I_{\{V_{u,v_i} \}}$,  where $M_n(u)$ is a random variable. Note that $ \mathbb{E}[M_n(u)] = \prod_{i=jn+1}^{jn+n}d_j $. Taking expectation and using principle of substitution fi\-ni\-shes the proof. \qed

\noindent
\textit{Proof of Theorem~\ref{T: esferic}}\\
Assume that $D(\bbT_F) > 0$. Then, for all $\alpha \in
(0, D(\bbT_F))$ there exists $N = N(\alpha)$ such that for all $n \geq N$
\begin{displaymath}
\min_{i \in \mathbf{N}}\frac{1}{n}\ln \left [\prod_{j=i+1}^{i+n}d_j \right]  > \alpha
\end{displaymath}
where
\begin{equation} \label{E: DF}
\prod_{j=i+1}^{i+n}d_j \geq e^{\alpha n} \textrm { for all } i \in \mathbf{N}
\textrm { and } n \geq N.
\end{equation}

Now we write
\begin{displaymath}
Z_{j+1} = \sum_{i=1}^{Z_j}Y_{j,i}^n,
\end{displaymath}
where $Y_{j,i}^n$ are i.i.d. copies of $Y_j^n$, being the number of offspring from the $i-th$ individual of the $j-th$ generation. By considering
Lemma~\ref{L: aes1} we have for all $j$ that
\begin{displaymath}
\mathbb{E}\left[\frac{Y_j^n}{\mu_j}\right]  = \frac{1}{\rho_j^{(n)}} \leq (\bbP [R > 0])^{-n}
\end{displaymath}
where $\rho_{j}^{(n)} = \bbP (V_{u,v})$, for any fixed pair $u \leq v$ such that $d({\mathcal O},u) =jn$
and $d({\mathcal O},v)=(j+1)n.$

Besides, by \ref{E: DF},  Lemma~\ref{L: aes2} and Lemma~\ref{L: aes1}
\begin{displaymath}
\textrm {if } \lim_{n \rightarrow \infty}
\sqrt[n]{\rho_n} > e^{-\alpha} \textrm { then } \liminf_{j \rightarrow \infty} \mu_j > 1.
\end{displaymath}

So, from Theorem 1 in Souza $\&$ Biggins~(\cite[p.40]{SouzaBiggins}), we conclude that the cone percolation process has a giant component with positive
probability if
\begin{displaymath}
\lim_{n \rightarrow \infty}
\sqrt[n]{\rho_n} > e^{-\alpha}.
\end{displaymath}
As this holds for every $\alpha \in (0, D(\bbT_F))$, the condition
\begin{displaymath}
\lim_{n \rightarrow \infty}\sqrt[n]{\rho_n} > e^{-D(\bbT_F)}
\end{displaymath}
guarantees the survival of the process with positive probability. \qed

\subsubsection{Homogeneous Galton-Watson Branching Trees}
\hfill

\noindent
\textit{Proof of Theorem~\ref{T: HGWT}}\\
Let us define two auxiliary branching process, being the first one
$\{\mathcal{X}_n\}_{n \in \bbN}$. For this process,
\begin{align*}
\mathbb{E}(\mathcal{X}) = \sum_{n=0}^{\infty}\bbP(R = n)\mathbb{E}(\mathcal{X} | R = n)
\end{align*}
where
\begin{align*}
\mathbb{E}(\mathcal{X} | R = 0) &= 0 \\
\mathbb{E}(\mathcal{X} | R = n) &= d^n, \hbox{ for } n =1,2, \dots.
\end{align*}
Note that
\begin{align}
\label{E: Minorante}
\mathbb{E}(\mathcal{X}) = \mathbb{E}[d^R]-p_0.
\end{align}

The second auxiliary process is $\{\mathcal{Y}_n\}_{n \in \bbN}$.
For this process
\begin{align*}
\mathbb{E}(\mathcal{Y}) = \sum_{n=0}^{\infty}\bbP(R = n)\mathbb{E}(\mathcal{X} | R = n)
\end{align*}
where 
\begin{align*}
\mathbb{E}(\mathcal{Y} | R = n) = d + d^2 + \cdots d^n.
\end{align*}
Note that
\begin{align}
\label{E: Majorante}
\mathbb{E}(\mathcal{Y}) = \frac{d}{d-1}(\mathbb{E}[d^R]-1).
\end{align}

Firstly we can assure \textit{(I)} by a comparison with
a supercritical branching process. In order to prove
\textit{(II)} and the left hand side one can see that our process dominates
$\{\mathcal{X}_n\}_{n \in \bbN}$. This process survives as long as
$\mathbb{E}[X]> 1$ therefore from~(\ref{E: Minorante}) our process survives
if $\mathbb{E}[d^R] > 1 + p_0.$

Secondly, also by a coupling argument, our process is do\-mi\-na\-ted by
$\{\mathcal{Y}_n\}_{n \in \bbN}.$ That process dies out provided $\mathbb{E}[Y] \leq 1$
and $\bbP[Y=1] \not= 1$,
therefore from~(\ref{E: Majorante}) our process dies out if
$\mathbb{E}[d^R] \leq 2 - \frac{1}{d},$ proving \textit{(III)} and the right hand side.

The proof of \textit{(IV)}  follows from the fact that
\[ \frac{1}{1-\mathbb{E}[\mathcal{X}]} \leq \mathbb{E}[|I|] \leq \frac{1}{1-\mathbb{E}[\mathcal{Y}]}. \]

\noindent
Acknowledgments: F.P.M. wishes to thank NYU-Shanghai China
and V.V.J. and K.R. wish to thank Instituto de Matem\'atica e 
Estat\'{\i}stica-USP Brazil
for kind hospitality.

\end{document}